\documentclass[10pt]{article}
\usepackage{anysize}
\marginsize{2.0cm}{2.0cm}{2.0 cm}{2.0 cm}
\usepackage{setspace}

\usepackage{latexsym,amsfonts,amsmath}
\usepackage{graphicx}
\usepackage{color}
\usepackage{epstopdf}

\newtheorem{definition}{Definition}[section]

\newtheorem{lemma}{Lemma}[section]

\newtheorem{proposition}{Proposition}[section]
\newtheorem{theorem}{Theorem}[section]
\newtheorem{remark}{Remark}[section]
\newtheorem{example}{Example}[section]

\def\PP{\mathsf{P}}
\def\EE{\mathsf{E}}

\begin{document}

\title{Further remarks on absorbing Markov decision processes}

\author{  Yi
Zhang \thanks{Corresponding author. School of Mathematics, University of Birmingham, Edgbaston,
Birmingham,
B15 2TT, U.K.. Email: y.zhang.29@bham.ac.uk} \and Xinran Zheng \thanks{School of Mathematics, University of Birmingham, Edgbaston,
Birmingham,
B15 2TT, U.K. E-mail: xxz299@student.bham.ac.uk}}
\maketitle

\par\noindent\textbf{Abstract}. In this note, based on the recent remarkable results of Dufour and Prieto-Rumeau, we deduce that for an absorbing {Markov decision process} with a given initial state, under a standard compactness-continuity condition, the space of occupation measures has the same convergent sequences, when it is endowed with the weak topology and with the weak-strong topology.   We provided two examples demonstrating that imposed condition cannot be replaced with its popular alternative, and the above assertion does not hold for the space of marginals of occupation measures on the state space. Moreover, the examples also clarify some results in the previous literature.

\bigskip


\section{Introduction}
In this note we provide further observations and examples concerning absorbing Markov decision processes (MDPs) with Borel state and action spaces.

An MDP model is called absorbing (with a given initial state), if there is an isolated cemetery being singled out in the state space, and under each strategy, the expected hitting time to the cemetery is finite, see e.g., \cite{Altman:1999,FeinbergRothblum:2012}. The expected hitting time can be written as a series of tail probabilities of the hitting time. If this series converges uniformly over the set of all strategies, then an absorbing MDP model is called uniformly absorbing. This definition comes from \cite{FeinbergPiunovskiy:2019}, where it was demonstrated that an absorbing MDP model is not necessarily uniformly absorbing. Uniformly absorbing MDP models with Borel state and action spaces have been studied in \cite{FeinbergRothblum:2012,FeinbergPiunovskiy:2019}, and they possess additional desirable properties. One of these is that, under either of two standard compactness-continuity conditions, the space of occupation measures in a uniformly absorbing MDP model is compact in the space of finite Borel measures on the product of the state and action spaces with respect to the weak topology (abbreviated henceforth the w-topology). In general, this property can be lost if an absorbing MDP model is not uniformly absorbing, as was demonstrated in \cite{Dufour:2023,PiunovskiyZhang:2023}.  The aforementioned compactness-continuity conditions are referred to as (S) and (W), and they are widely used in the literature of MDPs, see e.g, \cite{Schal:1975}. Both of them require the action space to be compact, whereas Condition (S) requires the transition kernel to be setwise continuous with respect to the actions for each fixed state, and Condition (W) requires the transition kernel to be weakly continuous with respect to the state-action pairs, see Section \ref{ZhangZhengYear2SectionModel} for their precise formulations.  Condition (S) is neither  implied by nor implies Condition (W).

Recently, there have been interests in characterizing when an absorbing MDP model is uniformly absorbing.   In \cite{Dufour:2024}, it was shown that, under either Condition (S) or Condition (W) (in fact under a more general condition that can be satisfied by them both), the w-compactness of the space of occupation measures is actually necessary and sufficient for an absorbing MDP model to be uniformly absorbing. Another necessary and sufficient condition was obtained in \cite{PiunovskiyZhang:2023}, being that the projection mapping from the space of strategic measure to the space of occupation measures is continuous, when both the spaces are endowed with their w-topology. In \cite{Dufour:2023}, under Condition (S), {yet} another characterization was given, which is the compactness of the space of occupation measures with respect to the weak-strong topology (ws-topology). In general, the ws-topology is finer than the w-topology, both of which are defined in Section \ref{ZhangZhengYear2SectionModel}. We also mention that in \cite{Dufour:2023}, the state space was allowed to be an arbitrary measurable space without involving topological properties, and therefore, Condition (W) was irrelevant therein.

The present note provides further observations and examples concerning absorbing MDPs, under Condition (S). The contributions are as follows. First, we show that in an absorbing MDP model, a sequence of occupation measures converges to an occupation measure in the w-topology if and only if the convergence takes place in the ws-topology (see Proposition \ref{ZhangZheng2024year2Theorme001}), and we provide an example (namely, Example \ref{ZhangZheng2024Example01}) demonstrating that this result does not hold if Condition (S) is replaced with Condition (W). We also present a second example (namely, Example \ref{ZhangZheng2024Example02}) demonstrating that under Condition (S), the previous result does not hold for the marginals of the occupation measures on the state space.  

A second contribution of this note can be seen as follows. Recall that in the definition of uniformly absorbing MDPs, the set of all strategies is considered. If, in that definition one is restricted to a given subset of strategies, say $\Lambda$, then one arrives at, what we call, $\Lambda$-uniformly {absorbing} MDP models. The consideration of the resulting notion is needed when we establish the first result. It is based on the observation that, under Condition (S), if the space of occupation measures over $\Lambda$ is w-compact, then an absorbing MDP model is $\Lambda$-uniformly absorbing, see Lemma \ref{ZhangZheng2024Lemma01}. {It is worth mentioning that the opposite direction of the implication in the last assertion does not hold in general, see Remark \ref{ZhangZheng2024RevisionRemark01}.}  On the other hand,  it is also possible to draw from the reasoning in the proof of \cite[Theorem 4.10]{Dufour:2024} that, under Condition (S), an absorbing MDP model  is $\Lambda$-uniformly absorbing if the space of marginals of occupation measures on the state space is compact with respect to the topology of setwise convergence (the s-topology). It is natural to ask whether the w-compactness of  the space of marginals of occupation measures over $\Lambda$ is sufficient for the absorbing MDP model {to be} $\Lambda$-uniformly absorbing. Example \ref{ZhangZheng2024Example02} also demonstrates that the answer to this question is negative. Therein lies our second contribution. 

Finally,  an immediate consequence of Proposition \ref{ZhangZheng2024year2Theorme001} is that for uniformly absorbing MDPs, the ws-topology and the w-topology on the space of occupation measures for a given initial state are equal, where we mention that occupation measures are infinite sums of certain marginals of strategic measures. A related result to this one is \cite{Nowak:1988}, which shows that under Condition (S), the w-topology and ws${}^\infty$-topology coincide on the space of strategic measures.  Example \ref{ZhangZheng2024Example01} also demonstrates that Condition (S) cannot be replaced with Condition (W) in that result. {In Theorem \ref{ZhangZheng2024RevisionTheorem01}, we also provide a further extension of the  result mentioned in the beginning of this paragraph from a given initial state to a finite set of initial states, and demonstrate in Remark \ref{ZhangZheng2024RevisionRemark02} that this extension may fail to hold for a more general set of initial states.}

The rest of the paper is organized as follows. In Section \ref{ZhangZhengYear2SectionModel}, we describe the MDP model, introduce Conditions (S) and (W), and recall the definitions of the w-, ws-, s-topologies. In Section \ref{ZhangZheng2024Year2Sec02}, we present the main results including the two examples. This note is ended with concluding remarks in Section \ref{ZhangZhengYear2Conclusion}.

\section{Model description and definitions}\label{ZhangZhengYear2SectionModel}

By an MDP model we mean the following collection $\{{\bf X},{\bf A},p\}$ of system primitives, where
\begin{itemize}
\item $\bf X$ and $\bf A$ are the state and action spaces, both assumed to be nonempty topological Borel spaces, endowed with their Borel $\sigma$-algebras ${\cal B}({\bf X})$ and ${\cal B}({\bf A})$, respectively;
\item $p(dy|x,a)$ is the transition stochastic kernel, a (measurable) stochastic kernel on ${\cal B}({\bf X})$ given ${\bf X}\times{\bf A}$.
\end{itemize}
Here, recall that  a Borel space is a Borel subset of some complete metrizable separable space.

The rest of this section is about a given MDP model $\{{\bf X},{\bf A},p\}$.
We first introduce some notations in order to define strategies and to describe the construction of the controlled process $\{X_n\}_{n=0}^\infty$ and the action process $\{A_n\}_{n=1}^\infty$.
Let the space of trajectories (or say histories) be the following countably infinite product
 $
  {\bf H}:={\bf X}\times({\bf A}\times{\bf X})^\infty.
 $
The generic notation for an element of ${\bf H}$ is $h=(x_0,a_1,x_1,\ldots)$. We endow ${\bf H}$ with the product topology and the product $\sigma$-algebra, which is also the Borel $\sigma$-algebra on it.
Let \begin{eqnarray*}
X_n(h)=x_n, A_{n+1}(h)=a_{n+1},~\forall~n\ge 0.
\end{eqnarray*}
They are the state and action variables at stage $n\ge 0$.

Let ${\bf H}_0:={\bf X}$, and for each $n\ge 1,$
$
{\bf H}_n:=({\bf X}\times {\bf A})^n\times {\bf X},
$
{endowed with the product $\sigma$-algebra}.
The generic notation for an element of ${\bf H}_0$ is $h_0=x_0$.
For each $n\ge 1$, the generic notation for an element of ${\bf H}_n$ is $h_n=(x_0,a_1,\dots,x_{n-1},a_{n},x_n)$.

\begin{definition}
If for each $n\ge 1,$ $\pi_n(da|h_{n-1})$ is a stochastic kernel on ${\cal B}({\bf A})$ given ${\bf H}_{n-1}$, then the sequence
 $\pi=\{\pi_n\}_{n=1}^{\infty}$ is called a strategy. {A strategy $\pi=\{\pi_n\}_{n=1}^{\infty}$ is called deterministic Markov if for each $n\ge 1$, there is a measurable mapping $\psi_n$ from ${\bf X}$ to ${\bf A}$ such that $\pi_n(da|h_{n-1})=\delta_{\psi_{n}(x_{n-1})}(da)$ for all $x_{n-1}\in {\bf X}$, where $\delta_{\psi_{n}(x_{n-1})}(da)$ is the Dirac measure on $({\bf A},{\cal B}({\bf A}))$ concentrated on the singleton $\{\psi_n(x_{n-1})\}$. 
 We identify such a deterministic Markov strategy with the underlying sequence of measurable mappings $\{\psi_n\}_{n=1}^\infty$. A deterministic Markov strategy $\pi=\{\psi_n\}_{n=1}^\infty$ is called deterministic stationary if for some measurable mapping $\psi$ from ${\bf X}$ to ${\bf A}$, $\psi_n=\psi$ for all $n\ge 1$.}  We identify such a deterministic stationary strategy with the underlying measurable mapping $\psi$.
\end{definition}
 The {sets} of all strategies {and all deterministic Markov strategies are} denoted by $\Pi$ {and $\Pi_{DM}$, respectively}.

In what follows, let the initial state $x_0\in  {\bf X}$ be given, {unless stated otherwise}. If a strategy $\pi$ is also fixed, then the strategic measure on $({\bf H},{\cal B}(\bf H))$, constructed in the standard way using the Ionescu-Tulcea theorem, is denoted as $\PP^\pi_{x_0}$. It is the unique probability measure on $({\bf H},{\cal B}({\bf H}))$ such that
\begin{eqnarray*}
&&\PP^\pi_{x_0}(X_0\in dx)=\delta_{x_0}(dx);\\
&& \PP^\pi_{x_0}(A_{n+1}\in da|X_0,A_1,\dots,X_n)=\pi_{n+1}(da|X_0,A_1,\dots,X_n);\\
&& \PP^\pi_{x_0}(X_{n+1}\in dx|X_0,A_1,\dots,X_n,A_{n+1})=p(dx|X_n,A_{n+1})~\forall~n\ge 0.
\end{eqnarray*}
The corresponding expectation is denoted as $\EE_{x_0}^\pi$.  It is convenient to introduce the following notation:
\begin{eqnarray*}
\PP_{x_0}^\pi|_{{\bf H}_n}(dh_n):=\PP_{x_0}^\pi((X_0,A_1,\dots,X_{n})\in dh_n)~\forall~n\ge 0.
\end{eqnarray*}

Let
\begin{eqnarray*}
{\cal P}_{x_0}:=\{\PP_{x_0}^\pi:\pi\in \Pi\}
\end{eqnarray*}
be the space of all strategic measures with the given initial state $x_0.$ Here we consider the case of an initial state instead of a more general initial distribution $\PP_0$ for the matter of brevity.

Now we present two standard compactness-continuity conditions that are  widely used in the literature of MDPs, see e.g., \cite{Schal:1975}.
\bigskip

\par\noindent\textbf{Condition (W)} \textit{ The action space ${\bf A}$ is compact, and for each bounded continuous function $g$ on ${\bf X}$, $\int_{{\bf X}}g(y)p(dy|x,a)$ is continuous in $(x,a)\in {\bf X}\times {\bf A}$.}

\bigskip

\par\noindent\textbf{Condition (S)} \textit{ The action space ${\bf A}$ is compact, and for each bounded measurable function $g$ on ${\bf X}$ and each $x\in {\bf X}$, $\int_{{\bf X}}g(y)p(dy|x,a)$ is continuous in $a\in  {\bf A}$.}
 \bigskip

We mention that in general, Condition (S) neither implies nor is implied by Condition (W).

Next, we describe a special class of MDP models $\{{\bf X},{\bf A},p\}$, where
${\bf X}={\bf Y}\cup\{\Delta\}$ with ${\bf Y}$ being a Borel space and $\Delta$ being an isolated point (cemetery), and
$
p(\{\Delta\}|\Delta,a)\equiv 1.
$
We remark that, as it stands, ${\bf Y}$ is allowed to contain other absorbing states, but $\Delta$ is singled out in our consideration throughout this paper.

For the MDP model $\{{\bf Y}\cup\{\Delta\},{\bf A},p\}$, let
$
\tau=\inf\{n\ge 0:~X_n=\Delta\}
$
be the hitting time of the state process $\{X_n\}_{n=0}^\infty$ {at the isolated cemetery $\Delta$}.

\begin{definition}
Consider the MDP model $\{{\bf Y}\cup\{\Delta\},{\bf A},p\}$, and let $\Lambda\subseteq \Pi$ be a nonempty set of strategies. The MDP model is called  absorbing with the initial state $x_0$ if, for each strategy $\pi\in\Pi$,
$
\EE^\pi_{x_0}[\tau]<\infty.
$
An absorbing MDP model with the initial state $x_0$ is called $\Lambda$-uniformly absorbing with the initial state $x_0$ if
$\lim_{n\to\infty}\sup_{\pi\in\Lambda} \EE^\pi_{x_0}\left[\sum_{t=n}^\infty \mathbb{I}\{\tau>t\}\right]=0,
$
where and below, $\mathbb{I}$ stands for the indicator function.
A $\Pi$-uniformly absorbing MDP model with the initial state $x_0$ is simply called uniformly absorbing.
\end{definition}

\begin{definition}
For the MDP model $\{{\bf Y}\cup\{\Delta\},{\bf A},p\}$, the occupation measure $\eta^\pi_{x_0}$ of a strategy $\pi$ with the initial state $x_0 $ is a measure   on $({\bf Y}\times {\bf A},{\cal B}({\bf Y}\times {\bf A}))
$ defined by
\begin{eqnarray*}
\eta_{x_0}^\pi(dx\times da):=\EE_{x_0}^\pi\left[\sum_{n=0}^\infty \mathbb{I}\{X_n\in dx, A_{n+1}\in da\}\right].
\end{eqnarray*}
\end{definition}
For the MDP model $\{{\bf Y}\cup\{\Delta\},{\bf A},p\}$ with the initial state $x_0$, for each $\Lambda\subseteq \Pi,$ the space of occupation measures over $\Lambda$  is denoted by
\begin{eqnarray*}
{\cal D}_{x_0}^\Lambda:=\left\{\eta_{x_0}^\pi:\pi\in\Lambda\right\}.
\end{eqnarray*}
Let
\begin{eqnarray*}
{{\cal D}_{x_0}^\Lambda}|_{\bf Y}:= \left\{{\eta_{x_0}^\pi}|_{\bf Y}:\pi\in\Lambda\right\},
\end{eqnarray*}
where for each measure $\eta$ on $({\bf Y}\times {\bf A},{\cal B}({\bf Y}\times {\bf A}))$, $\eta|_{\bf Y}(dx):=\eta(dx\times {\bf A})$ denotes its marginal on ${\bf Y}.$
 
For each $\pi\in\Pi$, $\eta_{x_0}^\pi({\bf Y}\times {\bf A})=\EE_{x_0}^\pi[\tau]$, so that if the MDP model $\{{\bf Y}\cup\{\Delta\},{\bf A},p\}$ is absorbing with the initial state $x_0$, then ${\cal D}_{x_0}^\Pi$  is a set of finite measures on $({\bf Y}\times {\bf A},{\cal B}({\bf Y}\times {\bf A})).$

Finally, we  recall some definitions to be used in the subsequent discussions.  Given a Borel space $\Omega$, let $\mathbb{P}(\Omega)$ be the space of Borel probability measures on $\Omega$, and ${\cal M}(\Omega)$ be the space of finite Borel measures on $\Omega.$ The w-topology on ${\cal M}(\Omega)$ is the coarsest topology on ${\cal M}(\Omega)$ with respect to which, for each bounded continuous function $g$ on $\Omega$, $\int_{\Omega}g(y)\eta(dy)$ is continuous in $\eta\in {\cal M}(\Omega).$ The w-topology on a subset of ${\cal M}(\Omega)$ is understood as the restriction of the w-topology on ${\cal M}(\Omega)$ to the subset.   According to \cite[Theorem 8.9.4(i)]{Bogachev:2007} and \cite[p.49]{Bogachev:2018}, the w-topology on  ${\cal M}(\Omega)$ is metrizable.

Now we turn to the MDP model $\{{\bf Y}\cup\{\Delta\},{\bf A},p\}$.
We call a real-valued measurable function $g$ on ${\bf Y}\times {\bf A}$ a Caratheodory function on ${\bf Y}\times {\bf A}$ if $g(x,a)$ is continuous in $a\in {\bf A}$ for each given $x\in{\bf Y}.$
 The ws-topology on ${\cal M}({\bf Y}\times {\bf A})$
  is the coarsest topology with respect to which, for each bounded Caratheodory function $g$ on ${\bf Y}\times {\bf A}$, $\int_{{\bf Y}\times {\bf A}}g(x,a)\eta(dx\times da)$ is continuous in $\eta\in {\cal M}({\bf Y}\times {\bf A}).$ According to \cite[p.496]{Balder:2001}, the ws-topology on ${\cal M}({\bf Y}\times {\bf A})$ is Hausdorff.

The s-topology on ${\cal M}({\bf Y})$ is  defined as the coarsest topology with respect to which, for each bounded measurable function $f$ on ${\bf Y}$, $\int_{\bf Y} f(x)\eta(dx)$ is continuous in $\eta\in {\cal M}({\bf Y})$, see \cite[Remark 1.21]{Ganssler:1971}.

\section{Main results}\label{ZhangZheng2024Year2Sec02}

\subsection{Statements}

In \cite[Theorem 3.1]{Dufour:2024}, it was shown that  if the MDP model $\{{\bf Y}\cup\{\Delta\}, {\bf A},p\}$ is absorbing with the initial state $x_0$ and satisfies Condition (S), then it is uniformly absorbing if ${\cal D}_{x_0}^\Pi$ is w-compact. The reasoning in the proof therein can be applied to deducing the following result.

\begin{lemma}\label{ZhangZheng2024Lemma01}
Suppose that the MDP model $\{{\bf Y}\cup\{\Delta\}, {\bf A},p\}$ is absorbing with the initial state $x_0 $, and satisfies Condition (S). Let $\Lambda$ be a nonempty set of strategies. Then the MDP model is uniformly absorbing over $\Lambda$ if ${\cal D}_{x_0}^\Lambda$ is  compact in ${\cal M}({\bf Y}\times {\bf A})$ with respect to the w-topology.
\end{lemma}

\begin{remark}\label{ZhangZheng2024RevisionRemark01} 
{According to \cite[Theorem 3.1]{Dufour:2024}, when $\Lambda=\Pi$, the opposite direction of the implication in the statement of Lemma \ref{ZhangZheng2024Lemma01} holds, too, i.e., under Condition (S),  if the MDP model is uniformly absorbing with the initial state $x_0$, then ${\cal D}_{x_0}^\Pi$ is  compact in ${\cal M}({\bf Y}\times {\bf A})$ with respect to the w-topology. However, this assertion may fail to hold for proper subsets $\Lambda$ of $\Pi$. As an example, consider the MDP model  $\{{\bf Y}\cup\{\Delta\}, {\bf A},p\}$, where ${\bf Y}=[0,1]\cup\{x_0\}$, where the relative topology on $[0,1]$ is Euclidean, and $x_0\notin [0,1]$ is an isolated point in ${\bf Y}$, ${\bf A}=\{0,1\}$ is endowed with the discrete topology, $p(\Gamma|x_0,a)\equiv Leb(\Gamma)$ for each $\Gamma\in {\cal B}([0,1])$ with $Leb(dy)$ being the Lebesgue measure on $([0,1],{\cal B}([0,1]))$, and $p(\{\Delta\}|x,a)=1$ for all $x\in [0,1]$ and $a\in\{0,1\}$. This MDP model satisfies Condition (S) and is uniformly absorbing (and thus $\Pi_{DM}$-uniformly absorbing), whereas from the reasoning in \cite[p.170-171]{Piunovskiy:1997}, one can see that ${\cal D}_{x_0}^{\Pi_{DM}}$ is not compact in ${\cal M}({\bf Y}\times {\bf A})$ in the w-topology.}
\end{remark}

\begin{lemma}\label{MorningZhangZheng2024Lemma01}
Suppose that the MDP model $\{{\bf Y}\cup\{\Delta\}, {\bf A},p\}$ is absorbing with the initial state $x_0\in{\bf Y}$, and satisfies Condition (S). Let $\Lambda$ be a nonempty set of strategies.  Then the MDP model is $\Lambda$-uniformly absorbing with the initial state $x_0$ if and only if ${{\cal D}_{x_0}^\Lambda}$ is relatively compact in ${\cal M}({\bf Y}\times {\bf A})$ with respect to the ws-topology.
\end{lemma}
\par\noindent\textit{Proof.} The reasoning in the proof of (a)$\Rightarrow$ (b) in \cite[Theorem 4.10]{Dufour:2023} proves the necessity part. {For the sufficiency part, assume that $\mathcal{D}^\Lambda_{x_0}$ is relatively ws-compact in ${\cal M}({\bf Y}\times {\bf A})$. By \cite[Lemma 4.9]{Dufour:2023}, ${\cal D}_{x_0}^\Pi$ is ws-closed in ${\cal M}({\bf Y}\times {\bf A})$. Take the closure of $\mathcal{D}^\Lambda_{x_0}$ in ${\cal M}({\bf Y}\times{\bf A})$ endowed with the ws-topology, say $cl({\mathcal{D}}^\Lambda_{x_0}).$ Since ${\cal D}_{x_0}^\Pi$ is ws-closed in ${\cal M}({\bf Y}\times {\bf A})$,  we see $cl({\mathcal{D}}^\Lambda_{x_0})\subseteq {\cal D}_{x_0}^\Pi$. Since  $\mathcal{D}^\Lambda_{x_0}$ is relatively ws-compact in ${\cal M}({\bf Y}\times {\bf A})$, we see that $cl({\mathcal{D}}^\Lambda_{x_0})$ is ws-compact and thus w-compact in ${\cal M}({\bf Y}\times {\bf A})$. Recall that the w-topology is coarser than the ws-topology on ${\cal M}({\bf Y}\times {\bf A})$. Let $\Lambda'$ be the collection of all strategies, whose occupation measures belong to $cl({\mathcal{D}}^\Lambda_{x_0})$. Then $\Lambda\subseteq\Lambda'$. By Lemma \ref{ZhangZheng2024Lemma01} applied to $\Lambda'$, we see that the MDP model is $\Lambda'$-uniformly absorbing, and thus $\Lambda$-uniformly absorbing, as required.} 
$\hfill\Box$


\begin{proposition}\label{ZhangZheng2024Theorem02}
Suppose that the MDP model $\{{\bf Y}\cup\{\Delta\}, {\bf A},p\}$ satisfies Condition (S), and is absorbing with the initial state $x_0$. Then a sequence $\{\eta_n\}_{n\ge 1}\subseteq {\cal D}_{x_0}^\Pi$ converges to some $\eta\in {\cal D}_{x_0}^\Pi$ in the w-topology if and only if $\eta_n\rightarrow \eta$ in the ws-topology.
\end{proposition}
\par\noindent\textit{Proof.} Consider a sequence $\{\eta_n\}_{n\ge 1}\subseteq {\cal D}_{x_0}^\Pi$ that converges to some $\eta\in {\cal D}_{x_0}^\Pi$ in the w-topology. We only need argue that $\{\eta_n\}_{n\ge 1}$ converges to $\eta$ in the ws-topology, because the ws-topology is finer than the w-topology.

Let $\pi_n$ ($n\ge 1$) and $\pi$ be strategies satisfying $\eta_n=\eta_{x_0}^{\pi_n}$ ($n\ge 1$) and $\eta=\eta_{x_0}^{\pi}$. Put $\Lambda=\{\pi_n,~n\ge 1\}\cup\{\pi\}$. Then ${\cal D}_{x_0}^\Lambda=\{\eta_n,~n\ge 1\}\cup\{\eta\}$ is compact in ${\cal M}({\bf Y}\times {\bf A})$ with respect to the w-topology. Since the w-topology on ${\cal M}({\bf Y}\times {\bf A})$ is metrizable, ${\cal D}_{x_0}^\Lambda$ is w-closed. It follows that ${\cal D}_{x_0}^\Lambda$ is ws-closed, because the ws-topology is finer than the w-topology.

By Lemma \ref{ZhangZheng2024Lemma01}, the MDP model is $\Lambda$-uniformly absorbing with the initial state $x_0$. By Lemma \ref{MorningZhangZheng2024Lemma01}, ${\cal D}_{x_0}^\Lambda$ is relatively ws-compact. This together with ${\cal D}_{x_0}^\Lambda$ being ws-closed  implies that ${\cal D}_{x_0}^\Lambda$ is ws-compact. According to \cite[Proposition 8.10.64]{Bogachev:2007}, it follows that ${\cal D}_{{x_0}}^\Lambda$ is metrizable, and thus sequentially ws-compact. Now any subsequence of $\{\eta_n\}_{n=1}^\infty\subseteq {\cal D}_{x_0}^\Lambda$ has a further ws-convergent subsequence, whose limit is necessarily $\eta$, because $ \eta_n\rightarrow \eta$ in the w-topology, which is coarser than then ws-topology. This shows that $ \eta_n\rightarrow \eta$ in the ws-topology, as required.
 $\hfill\Box$

{For each $F\subseteq {\bf X}$, introduce ${\cal D}_{F}^{\Pi}:=\{\eta\in{\cal D}_{x}^\pi:~ x\in F,~\pi\in\Pi\}$, the set of all occupation measures for all initial states from $F$.}
\begin{theorem}\label{ZhangZheng2024RevisionTheorem01}
{Suppose that the MDP model $\{{\bf Y}\cup\{\Delta\},{\bf A}, p\}$ with ${\bf X}={\bf Y}\cup\{\Delta\}$ satisfies Condition (S) and is uniformly absorbing for each initial state $x\in F$, where $F$ is a nonempty finite subset of ${\bf X}$, so that ${\cal D}_F^{\Pi}\subseteq {\cal M}({\bf Y}\times {\bf A})$. Then the w-topology and ws-topology on ${\cal D}_F^\Pi$ coincide.} 
\end{theorem}

\par\noindent\textit{Proof.} {According to \cite[Theorem 4.10]{Dufour:2023}, under Condition (S), if the MDP model is uniformly absorbing with the initial state $x\in F$, then ${\cal D}_{x}^\Pi$ is ws-compact. Since $F$ is finite, ${\cal D}^\Pi_F$ is ws-compact, as a finite union of compact sets. This implies, by \cite[Proposition 8.10.64]{Bogachev:2007}, that ${\cal D}_{F}^\Pi$ is metrizable. Since the w-topology on ${\cal D}_{F}^\Pi$ is metrizable, too, it suffices to prove that ${\cal D}_{F}^\Pi$ has the same convergent sequences in w-topology and in ws-topology, and this amounts to showing that if a sequence $\{\eta_n\}_{n\ge 1}\subseteq {\cal D}_{F}^\Pi$ converges to some $\eta\in {\cal D}_{F}^\Pi$ in the w-topology, then $\eta_n\rightarrow \eta$ in the ws-topology.  Assume for contradiction that this convergence does not take place in the ws-topology. Recall that a sequence converges to a point if each of its subsequences has a further  subsequence converging to that point. Then $\{\eta_n\}_{n\ge 1}\subseteq {\cal D}_{F}^\Pi$ has a subsequence $\{\eta_{n_k}\}_{k\ge 1}$ without a further subsequence ws-converging to $\eta$. Since $F$ is finite, for some $x_0\in F$, $\{\eta_{n_k}\}_{k\ge 1}$ has a subsequence $\{\eta_{n_{k_j}}\}_{j\ge 1} \subseteq{\cal D}_{x_0}^\Pi$. By assumption, $\{\eta_{n_{k_j}}\}_{j\ge 1}$ w-converges to $\eta$. Since the MDP model is uniformly absorbing with the initial state $x_0\in F$, by Remark \ref{ZhangZheng2024RevisionRemark01}, ${\cal D}_{x_0}^\Pi$ is w-compact, and hence $\eta\in {\cal D}_{x_0}^\Pi$. Now, Proposition \ref{ZhangZheng2024Theorem02} implies that $\{\eta_{n_{k_j}}\}_{j\ge 1}\subseteq {\cal D}_{x_0}^\Pi$ ws-converges to $\eta\in {\cal D}_{x_0}^\Pi$, which is a desired contradiction. $\hfill\Box$}

\begin{remark}\label{ZhangZheng2024RevisionRemark02} 
{Theorem \ref{ZhangZheng2024RevisionTheorem01} may fail to hold when $F$ is a countable  compact set. Consider the MDP model with ${\bf Y}=\{\frac{1}{n},~n\ge 1\}\cup\{0\}$, which is endowed with the relative Euclidean topology and compact,  ${\bf A}=\{a^\ast\}$ being a singleton, and $p(\{\Delta\}|x,a)\equiv 1$. It satisfies Condition (S) and is uniformly absorbing for each initial state. We have ${\cal D}^\Pi_{\bf Y}=\{\delta_{\frac{1}{n}}(dy)\delta_{a^\ast}(da),~n\ge 1\}\cup\{\delta_0(dy)\delta_{a^\ast}(da)\}.$ As $n\rightarrow \infty$, $\delta_{\frac{1}{n}}(dy)\delta_{a^\ast}(da)$ converges to $\delta_0(dy)\delta_{a^\ast}(da)$ in the w-topology but not in the ws-topology.}
\end{remark}

\subsection{Examples}
\subsubsection{On the role of Condition (S)}
The first example demonstrates that  the assertion in Proposition \ref{ZhangZheng2024Theorem02} may fail to hold, if Condition (S) is therein replaced with Condition (W).

\begin{example}\label{ZhangZheng2024Example01} The state space is given by ${\bf X}:={\bf Y}\cup\{\Delta\},$ where ${\bf Y}:=\{0,1\}\times [0,1]$. Here, $[0,1]$ is endowed with the Euclidean topology, $\{0,1\}$ is endowed with the discrete topology, and ${\bf Y}$ is endowed with the product topology. The action space is given by ${\bf A}:=[0,1]$, also endowed with the Euclidean topology. The transition kernel is fully specified for all $a\in {\bf A}$ by
\begin{eqnarray*}
&&p(\{1\}\times dy|(0,x),a):=\delta_a(dy) ~\mbox{on ${\cal B}([0,1])$}~\forall~x\in [0,1];\\
&&p(\{\Delta\}|(1,x),a):=1~\forall~x\in [0,1];\\
&&p(\{\Delta\}|\Delta,a):=1.
\end{eqnarray*}
The initial state is $x_0=(0,0)$.
\end{example}

\begin{proposition}\label{ZhangZheng2024year2Theorme001}
For the MDP model $\{{\bf Y}\cup\{\Delta\},{\bf A},p\}$ in Example \ref{ZhangZheng2024Example01} with the initial state $x_0=(0,0)$, the following assertions hold.
\begin{itemize}
\item[(a)] The MDP model satisfies Condition (W).
\item[(b)] The MDP model with the initial state $x_0$ is uniformly absorbing.
\item[(c)] There is a sequence $\{\eta_n\}_{n=1}^\infty\subseteq {\cal D}_{x_0}^\Pi$ that converges to some $\eta\in {\cal D}_{x_0}^\Pi$ in the w-topology but does not have any convergent subsequence in ${\cal M}({\bf Y}\times {\bf A})$ with respect to the ws-topology.
\item[(d)] The space ${\cal D}_{x_0}^\Pi$ is w-compact, but not ws-compact.
\end{itemize}
\end{proposition}

\par\noindent\textit{Proof.}  Consider the MDP model described in Example \ref{ZhangZheng2024Example01}. We use the generic notation $\bar{x}$ for an element of ${\bf X}={\bf Y}\cup\{\Delta\}$.

(a) This part holds trivially, because for each bounded continuous function $g$ on ${\bf Y}$,
\begin{eqnarray*}
\int_{{\bf Y}} g(y)p(dy|(n,x),a)=\mathbb{I}\{n=0\}g(1,a)~\forall~(n,x)\in {\bf Y}.
\end{eqnarray*}

(b) Since starting from $x_0=(0,0)$, the cemetery $\Delta$ is firstly reached by the state process exactly in two steps, under each strategy, we see that this MDP model is uniformly absorbing.

(c)  For each strategy $\pi=\{\pi_n\}_{n=1}^\infty$, with the initial state $x_0=(0,0)$, the occupation measure $\eta_{(0,0)}^\pi$ of $\pi$ is given by
\begin{eqnarray}\label{ZZYear2Eqn01}
&&\eta_{(0,0)}^\pi(\{0\}\times dx\times da)=\delta_{0}(dx)\pi_1(da|(0,0));\nonumber\\
&&
\eta_{(0,0)}^\pi(\{1\}\times dx\times da)=\pi_2(da|(0,0),x,(1,x))\pi_1(dx|(0,0)).
\end{eqnarray}

Now consider a sequence of strategies $\{\pi^m\}_{m=1}^\infty$, where for each $m\ge 1,$ $\pi^m=\{\pi_n^m\}_{n=1}^\infty$ is given for each $\bar{x}\in {\bf X}$ by
\begin{eqnarray*}
&&\pi^m_1(da|\bar{x})\equiv m\mathbb{I}\left\{a\in\left[0,\frac{1}{m}\right]\right\}da;\\
&&\pi^m_n(da|\bar{x})\equiv \mu(da)~ \forall~n\ge 2
\end{eqnarray*}
for some fixed $\mu\in\mathbb{P}({\bf A})$.
Let $\pi^\infty=\{\pi^\infty_n\}_{n=1}^\infty$ be a strategy given for each $\bar{x}\in {\bf X}$  by
\begin{eqnarray*}
&&\pi^\infty_1(da|\bar{x})= \delta_0(da);\\
&&\pi^\infty_n(da|\bar{x})= \mu(da)~ \forall~n\ge 2.
\end{eqnarray*}
Then $\eta_{(0,0)}^{\pi^m}$ converges to $\eta^{\pi^\infty}_{(0,0)}$ in the w-topology. Indeed, for each bounded continuous function $g$ on ${\bf Y}\times {\bf A}$, upon applying (\ref{ZZYear2Eqn01}) to $\pi^m$ and $\pi^\infty$, we see that
\begin{eqnarray}\label{ZZYear2Eqn02}
&&\int_{{\bf Y}\times {\bf A}} g(\bar{x},a)\eta_{(0,0)}^{\pi^m}(d\bar{x}\times da)\nonumber\\
&=&\int_{{\bf A}}g((0,0),a)\pi^m_1(da|(0,0))+\int_{[0,1]\times{\bf A}}g((1,x),a)\pi^m_2(da|(0,0),x,(1,x))\pi^m_1(dx|(0,0))\nonumber\\
&=&m\int_0^{\frac{1}{m}} g((0,0),a)da+m\int_0^{\frac{1}{m}} \left\{\int_{\bf A}g((1,x),a) \mu(da)\right\} dx,
\end{eqnarray}
which, as $m\rightarrow \infty,$ converges to
\begin{eqnarray}\label{ZZYear2Eqn05}
&&g((0,0),0)+ \left\{\int_{\bf A}g((1,0),a) \mu(da)\right\}\nonumber\\
&=&\int_{[0,1]} g((0,0),a)\delta_0(da)+\int_{[0,1]} \left\{\int_{\bf A}g((1,x),a) \mu(da)\right\} \delta_0(dx)\nonumber\\
&=&\int_{{\bf Y}\times {\bf A}} g(\bar{x},a)\eta_{(0,0)}^{\pi^\infty}(d\bar{x}\times da)
\end{eqnarray}
by the fundamental theorem of calculus, having observed that $a\in[0,1]\rightarrow g((0,0),a)$ and $x\in[0,1]\rightarrow \int_{\bf A}g((1,x),a) \mu(da)$ are bounded continuous.

Now we verify that $\{\eta_{x_0}^{\pi^m}\}_{m=1}^\infty$ does not have any convergent subsequence in ${\cal M}({\bf Y}\times {\bf A})$ with respect to the ws-topology. Suppose for contradiction that $\{\eta_{x_0}^{\pi^m}\}_{m=1}^\infty$ has a convergent subsequence in ${\cal M}({\bf Y}\times {\bf A})$ with respect to the ws-topology, say $\{\eta_{x_0}^{\pi^{m_k}}\}_{k=1}^\infty$.
Since the ws-topology is finer than the w-topology, and  $\{\eta_{x_0}^{\pi^m}\}_{m=1}^\infty$  converges to $\eta_{x_0}^{\pi^\infty}$ in the w-topology, we see that the limit of $\{\eta_{x_0}^{\pi^{m_k}}\}_{k=1}^\infty$, assumed to exist, is necessarily $\eta_{x_0}^{\pi^{\infty}}$. The rest verifies that $\{\eta_{x_0}^{\pi^{m_k}}\}_{k=1}^\infty$ does not converge to $\eta_{x_0}^{\pi^{\infty}}$ in the ws-topology.

Consider a function $g$ on ${\bf Y}\times {\bf A}$ defined by
\begin{eqnarray*}
g((n,x),a)=\hat{g}(x):=\mathbb{I}\{x>0\},~\forall~(n,x)\in{\bf Y}, ~a\in {\bf A}=[0,1].
\end{eqnarray*}
Then this is a bounded Caratheodory function on ${\bf Y}\times {\bf A}$, and is not continuous thereon.
Applying (\ref{ZZYear2Eqn02}) to this function $g$ and $\pi^{m_k}$ leads to
\begin{eqnarray*}
&&\int_{{\bf Y}\times {\bf A}} g(\bar{x},a)\eta_{(0,0)}^{\pi^{m_k}}(d\bar{x}\times da)\nonumber\\
&=&m_k\int_0^{\frac{1}{m_k}} g((0,0),a)da+m_k\int_0^{\frac{1}{m_k}} \left\{\int_{\bf A}g((1,x),a) \mu(da)\right\} dx\\
&=&m_k\int_0^{\frac{1}{m_k}}     \mathbb{I}\{x>0\}  dx=1
\end{eqnarray*}
for each $m_k\ge 1$. On the other hand, from (\ref{ZZYear2Eqn05}) applied to this function $g$, we see that
\begin{eqnarray*}
&&\int_{{\bf Y}\times {\bf A}} g(\bar{x},a)\eta_{(0,0)}^{\pi^\infty}(d\bar{x}\times da)\\
&=&g((0,0),0)+  \int_{\bf A}g((1,0),a) \mu(da) \nonumber\\
&=&2\hat{g}(0)=0.
\end{eqnarray*}
Thus, $\left\{\int_{{\bf Y}\times {\bf A}} g(\bar{x})\eta_{(0,0)}^{\pi^{m_k}}(d\bar{x}\times da)\right\}_{k=1}^\infty$ converges to $1 \ne 0=\int_{{\bf Y}\times {\bf A}} g(\bar{x})\eta_{(0,0)}^{\pi^\infty}(d\bar{x}\times da)$, and $\{\eta_{(0,0)}^{\pi^{m_k}}\}_{k=1}^\infty$ does not converge to $\eta_{(0,0)}^{\pi^\infty}$ in the ws-topology as $k\rightarrow \infty$. This is the required contradiction.

(d) The w-compactness of ${\cal D}_{x_0}^\Pi$ follows from \cite[Theorem 3.1]{Dufour:2024} or \cite[Theorem 2]{PiunovskiyZhang:2023} and \cite[Theorem 5.6]{Schal:1975}. If ${\cal D}_{x_0}^\Pi$ was ws-compact, then according to \cite[Proposition 8.10.64]{Bogachev:2007}, it would be metrizable and thus sequentially ws-compact, but this would then contradict part (c). $\hfill\Box$

Proposition \ref{ZhangZheng2024year2Theorme001}(c) shows that Condition (S) cannot be replaced with Condition (W) in Proposition \ref{ZhangZheng2024Theorem02}, and Proposition  \ref{ZhangZheng2024year2Theorme001}(d) shows that for a uniformly absorbing MDP {model}, the space of occupation measures may fail to be ws-compact under Condition (W), on contrary with the case under Condition (S), see \cite[Theorem 4.10]{Dufour:2023}.

\subsubsection{On Lemma \ref{ZhangZheng2024Lemma01} and Proposition \ref{ZhangZheng2024Theorem02}}

Let us observe that from Lemma \ref{MorningZhangZheng2024Lemma01}, one can deduce that, under Condition (S), an absorbing MDP model is $\Lambda$-uniformly absorbing with the initial state $x_0$ if   ${{\cal D}_{x_0}^\Lambda}|_{\bf Y}$ is compact in ${\cal M}({\bf Y})$ with respect to the s-topology. Indeed,  if ${{\cal D}_{x_0}^\Lambda}|_{\bf Y}$ is s-compact in ${\cal M}({\bf Y})$, then, according to \cite[Exercise 4.7.148(ii)]{Bogachev:2007}, it is metrizable, and thus sequentially s-compact in ${\cal M}(\bf Y)$.  By \cite[Theorems 2.5, 5.2]{Balder:2001}, it follows that ${{\cal D}_{x_0}^\Lambda}$ is relatively ws-compact. Hence, by Lemma \ref{MorningZhangZheng2024Lemma01}, the MDP model is $\Lambda$-uniformly absorbing with the initial state $x_0$.

In view of this fact, it is natural to ask whether in Lemma \ref{ZhangZheng2024Lemma01}, the w-compactness of ${\cal D}_{x_0}^\Lambda$ in ${\cal M}({\bf Y}\times {\bf A})$ can be replaced with the w-compactness of ${{\cal D}_{x_0}^\Lambda}|_{\bf Y}$ in ${\cal M}({\bf Y})$.

The next example demonstrates that this cannot be done in general. It also shows that the assertion of  Proposition \ref{ZhangZheng2024Theorem02} may not hold when the occupation measures $\eta_n,\eta$ are replaced with their marginals $\eta_n|_{{\bf Y}},\eta|_{{\bf Y}}$ on ${\bf Y}$: that $\eta_n|_{{\bf Y}}$ converges to $\eta|_{{\bf Y}}$ in the w-topology does not imply the convergence holds in the s-topology.

\begin{example}\label{ZhangZheng2024Example02}
The state space  is
given by $
  {\bf X}  ={\bf Y}\cup\{\Delta\}
$
with
$
{\bf Y}=\{b_n,~n=1,2,\dots\}\cup\{1\},
$
where
\begin{eqnarray*}
b_n:=\sum_{i=1}^{n}\frac{1}{2^i},~\forall~n\ge 1.
\end{eqnarray*}
We emphasize that here we endow ${\bf Y}$ with the (relative) Euclidean topology.
The action space  is
$
    {\bf A}=\{1,2,3\},
$
equipped with the discrete topology.
The transition kernel $p$ is defined as follows:
\begin{eqnarray*}
&&p(\{\Delta\}|\Delta,a)\equiv 1\equiv p(\{\Delta\}|1,a);\\
&&p(\{\Delta\}|b_1,1)=1,\\
&&p(\{b_n\}|b_n,1)=1-\frac{1}{2^{n-2}}, ~p(\{\Delta\}|b_n,1)=\frac{1}{2^{n-2}}~\forall~n\ge 2;
\end{eqnarray*}
and for all $n\ge 1,$
\begin{eqnarray*}
&&p(\{\Delta\}|b_n,2)=\frac{1}{2}=p(\{b_{n+1}\}|b_n,2);\\
&&p(\{\Delta\}|b_n,3)=\frac{1}{4}=p(\{1\}|b_n,3),~p(\{b_{n+1}\}|b_n,3)=\frac{1}{2}.
\end{eqnarray*}
We fix the initial state to be $x_0=b_1=\frac{1}{2}$.
Consider the following deterministic stationary strategies: for each $n\ge 3$, the deterministic stationary strategy $\psi^n$  is defined by
\begin{eqnarray*}
    \psi^n(x):=2\cdot\mathbb{I}\{x\in \{b_1,\dots,b_n\}\}+\mathbb{I}\{x\in\{b_{m},~m\ge n+1\}\}{+ 3\cdot\mathbb{I}\{x\in\{1,\Delta\}\}}~\forall~x\in {\bf Y}.
\end{eqnarray*}
Also, consider the deterministic stationary strategy $\psi$ defined by
$
\psi(x)\equiv 3.
$
Let \begin{eqnarray*}
\Lambda=\{\psi^n,~n\ge 3\}\cup\{\psi\}.
\end{eqnarray*}
 \end{example}

 In the above example,  if $a=1$ is chosen at $x=b_n$ with $n\ge 2$, in the next stage, the system stays at $b_n$ with probability $1-\frac{1}{2^{n-2}}$, and moves to the cemetery $\Delta$ with probability $\frac{1}{2^{n-2}}$. If $a=1$ is chosen at $x=b_1$, then, in the next stage, the system goes to $\Delta$ with probability $1$.
If $a=2$ is chosen at $x=b_n$ with $n\ge 1$, in the next stage, the system goes to  $b_{n+1}$ with probability $\frac{1}{2}$, and moves to $\Delta$ with probability $\frac{1}{2}$.
If $a=3$ is chosen at $x= b_n$ with $n\ge 1,$ then, in the next stage, the system goes to $b_{n+1}$ with probability $\frac{1}{2}$,  goes to state $1$ with probability $\frac{1}{4}$, and moves to $\Delta$ with probability $\frac{1}{4}$.
Any action chosen at state $1$ moves the system to $\Delta$ in the next stage with probability $1$.

\begin{proposition}
For the MDP model in Example \ref{ZhangZheng2024Example02} with the initial state $x_0=b_1=\frac{1}{2}$ and the set $\Lambda$ defined therein, the following assertions hold.
\begin{itemize}
\item[(a)] Condition (S) is satisfied.

\item[(b)]  The MDP model is absorbing but not $\Lambda$-uniformly absorbing with the initial state $x_0.$
\item[(c)] ${\eta^{\psi^n}_{x_0}}|_{\bf Y}$ converges to ${\eta^{\psi}_{x_0}}|_{\bf Y}$ in the w-topology, but not in the s-topology.
\item[(d)] ${{\cal D}_{x_0}^\Lambda}|_{\bf Y}$ is w-compact but not s-compact in ${\cal M}({\bf Y}).$
\end{itemize}
\end{proposition}

\par\noindent\textit{Proof.} Part (a) holds trivially because ${\bf A}$ was endowed with the discrete topology.

(b) Let us show that the MDP model $\{{\bf Y}\cup \{\Delta\},{\bf A},p\}$ is absorbing with the initial state $x_0=\frac{1}{2}.$   To this end, it suffices to show that the function $w^\ast$  defined on ${\bf X}$ is finite-valued, where
\begin{equation*}
    w^*(x):=\sup_{\pi\in\Pi}\EE^\pi_x \left[\sum^\infty_{t=0}\mathbb{I}\{X_t\in {\bf Y}\}\right]=\sup_{\pi\in \Pi} \eta^\pi_x( {\bf Y}\times{\bf A} )~{\forall~x\in {\bf X}={\bf Y}\cup\{\Delta\}}.
\end{equation*}
This is done as follows.

According to  \cite[Propositions 9.8 and 9.10]{Bertsekas:1978}, $w^\ast$ is the minimal nonnegative solution for the following Bellman equation
\begin{eqnarray}\label{ZhangZheng2024ExampleSeparateFileEqn01}
&&    w(\Delta)=0,\nonumber\\
&&    w(b_n)=1+\max\left\{\mathbb{I}\{n\geq 2\}\left(1-\frac{1}{2^{n-2}}\right)w(b_n),~ \frac{1}{2}w(b_{n+1}),~\frac{1}{2}w(b_{n+1})+\frac{1}{4}w(1)\right\}~\forall~n\ge 1,\nonumber\\
&&    w(1)=1,
\end{eqnarray}
where the three terms in the above parenthesis correspond to the action $a=1,2,3,$ respectively.
Let us verify that the $[0,\infty)$-valued function $w$ on ${\bf X}$ defined by $w(\Delta):=0$, $w(1):=1$ and
\begin{equation*}
    w(b_n):=\frac{1}{2}+2+2^{n-2}~\forall~ n\ge 1
\end{equation*}
is a solution for the Bellman equation (\ref{ZhangZheng2024ExampleSeparateFileEqn01}). Once this is verified, it would then follow that $w^\ast$ is finite-valued, because $w^\ast\le w$ pointwise. Now let us show that the function  $w$ satisfies the second equality in (\ref{ZhangZheng2024ExampleSeparateFileEqn01}) (the other equalities in (\ref{ZhangZheng2024ExampleSeparateFileEqn01}) are satisfied automatically). For this function $w$, the first term inside the parenthesis in (\ref{ZhangZheng2024ExampleSeparateFileEqn01}) can be computed as
\begin{eqnarray*}
    && \mathbb{I}\{n\geq 2\}\left(1-\frac{1}{2^{n-2}}\right)w(b_n)\\
    &=& \mathbb{I}\{n\geq 2\}\left(1-\frac{1}{2^{n-2}}\right)\left(\frac{1}{2}+2+2^{n-2}\right)\\
    &=& \mathbb{I}\{n\geq 2\}\left(\frac{1}{2}+1+2^{n-2}-\frac{1}{2^{n-1}}-\frac{1}{2^{n-3}}\right),
\end{eqnarray*}
the second term inside the parenthesis in (\ref{ZhangZheng2024ExampleSeparateFileEqn01}) reads
\begin{eqnarray*}
    &&\frac{1}{2}w(b_{n+1})=\frac{1}{2} \left(\frac{1}{2}+2+2^{n-1}\right)=\frac{1}{4}+1+2^{n-2},
\end{eqnarray*}
whereas  the third term inside the parenthesis in (\ref{ZhangZheng2024ExampleSeparateFileEqn01}) reads
\begin{eqnarray*}
    &&\frac{1}{2}w(b_{n+1})+\frac{1}{4}w(1)
    =\frac{1}{2}\left(\frac{1}{2}+2+2^{n-1}\right)+\frac{1}{4}=\frac{1}{2}+1+2^{n-2}.
\end{eqnarray*}
Hence, for the function $w$ under consideration, we have for each $n\ge 1$ that
\begin{eqnarray*}
&&    1+\max\left\{\mathbb{I}\{n\geq 2\}\left(1-\frac{1}{2^{n-2}}\right)w(b_n),~ \frac{1}{2}w(b_{n+1}),~\frac{1}{2}w(b_{n+1})+\frac{1}{4}w(1)\right\}\\
    &=&1+\max\left\{\mathbb{I}\{n\geq 2\}\left(\frac{1}{2}+1+2^{n-2}-\frac{1}{2^{n-1}}-\frac{1}{2^{n-3}}\right),~\frac{1}{4}+1+2^{n-2},~\frac{1}{2}+1+2^{n-2}\right\}\\
    &=&1+\frac{1}{2}+1+2^{n-2}=w(b_n),
\end{eqnarray*}
as required.

Next, we show that the model is not $\Lambda$-uniformly absorbing with the initial state $x_0=b_1=\frac{1}{2}$. Recall the definitions of $\Lambda$, $\psi^n$ and $\psi$ in Example \ref{ZhangZheng2024Example02}.

Observe that for each $n\ge 3$ and $x\in{\bf Y}$, ${\eta^{\psi^n}_{\frac{1}{2}}}|_{\bf Y}(\{x\})$ is given by
\begin{align}\label{ZhangZheng2024ExampleSeparateFileEqn02}
    {\eta^{\psi^n}_{\frac{1}{2}}}|_{{\bf Y}}(\{b_m\})=
    \begin{cases}
        \frac{1}{2^{m-1}}\quad\text{if}\quad   1\le m<n+1;\\\
        \frac{1}{2}\quad\text{if}\quad m=n+1;\\
        0\quad\text{otherwise}
    \end{cases}
\end{align}
and
\begin{eqnarray}\label{ZhangZheng2024ExampleSeparateFileEqn05}
{\eta^{\psi^n}_{\frac{1}{2}}}|_{{\bf Y}}(\{1\})=0.
\end{eqnarray}
Indeed, under the strategy $\psi^n$, for each $1\le m\le n$, state $b_m$ can be visited at most once, because the action $2$ is always in use at each of these states; and the probability of the system hitting $b_m$ is given by $\frac{1}{2^{m-1}}$, as the initial state is $x_0=b_1=\frac{1}{2}$. This justifies the first equality in (\ref{ZhangZheng2024ExampleSeparateFileEqn02}). For each $m>n+1$, $b_m$ is never visited, as well as state $1$. This justifies the last equality in (\ref{ZhangZheng2024ExampleSeparateFileEqn02}) as well as (\ref{ZhangZheng2024ExampleSeparateFileEqn05}). The probability of hitting state $b_{n+1}$ is $\frac{1}{2^n}$, and given that $b_{n+1}$ is hit, the duration of the system staying at $b_{n+1}$ obeys the geometric distribution with parameter $\frac{1}{2^{n-1}}$. Consequently,
\begin{eqnarray*}
   { \eta^{\psi^n}_{\frac{1}{2}}}|_{{\bf Y}}(\{b_{n+1}\})=
    \frac{1}{2^n} 2^{n-1}=\frac{1}{2},
\end{eqnarray*}
justifying the second equality in (\ref{ZhangZheng2024ExampleSeparateFileEqn02}).

Note  for all $n\ge 3$ that
\begin{equation*}
    \sup_{\pi\in \Lambda}\EE^\pi_{\frac{1}{2}}\left[\sum^\infty_{t=n}\mathbb{I}\{X_t\in {\bf Y}\}\right]\geq \EE^{\psi^n}_{\frac{1}{2}}\left[\sum^\infty_{t=n}\mathbb{I}\{X_t=b_{n+1}\}\right]={\eta^{\psi^n}_{\frac{1}{2}}}|_{\bf Y}(\{b_{n+1}\})=\frac{1}{2},
\end{equation*}
where the second to the last equality holds because, starting with the initial state $x_0=b_1=\frac{1}{2},$ under the strategy $\psi^n,$ the first time when $b_{n+1}$ is hit, if at all, can only be $n$.
Therefore, $\sup_{\pi\in\Lambda}\EE^\pi_{\frac{1}{2}}\left[\sum^\infty_{t=n}\mathbb{I}\{X_t\in {\bf Y}\}\right]$ does not converge to $0$ as $n\rightarrow\infty$, and the MDP model is not $\Lambda$-uniformly absorbing   with the initial state $x_0=b_1=\frac{1}{2}$.

(c,d) To show that  ${\cal D}_{x_0}^\Lambda$ is w-compact, it is sufficient to verify that ${\eta_{\frac{1}{2}}^{\psi^n}}|_{{\bf Y}}$ converges to ${\eta_{\frac{1}{2}}^{\psi}}|_{{\bf Y}}$  in the w-topology.

Direct calculations yield that
\begin{eqnarray*}
    &&{\eta^{\psi}_{\frac{1}{2}}}|_{{\bf Y}}(\{1\})=\EE^{\psi}_{\frac{1}{2}}\left[\sum^\infty_{t=0}\mathbb{I}\{X_t=1\}\right]=\sum^\infty_{t=1}\PP^{\psi}_{\frac{1}{2}}(X_t=1)\\
    &=& \PP_{\frac{1}{2}}^\psi(X_1=1)+\PP_{\frac{1}{2}}^\psi(X_1\ne \Delta, 1,~X_2=1)+\PP_{\frac{1}{2}}^\psi(X_1,X_2\ne \Delta, 1,~X_3=1) +\dots   \\
    &=&   \frac{1}{4}+\frac{1}{2}\cdot\frac{1}{4}+\frac{1}{2^2}\cdot\frac{1}{4}+\dots=\frac{1}{2},
\end{eqnarray*}
where the second equality holds because the initial state $x_0\ne 1$.
Moreover, similar considerations to those for ${\eta_{\frac{1}{2}}^{\psi^n}}|_{{\bf Y}}$ lead to
\begin{eqnarray}\label{ZhangZheng2024ExampleSeparateFileEqn06}
{\eta_{\frac{1}{2}}^\psi}|_{\bf Y}(\{b_n\})=\frac{1}{2^{n-1}} ~\forall~n\ge 1.
\end{eqnarray}
(Under the strategy $\psi,$ each $b_n$ with $n\ge 1$ can be visited at most once, and starting with $x_0=b_1=\frac{1}{2},$ the probability of hitting $b_n$ is $\frac{1}{2^{n-1}}$.)

Now let $f$ be an arbitrarily fixed bounded continuous function $f$ on ${\bf Y}$.
Then
\begin{eqnarray}\label{ZhangZheng2024ExampleSeparateFileEqn07}
    &&\left|\int_{{\bf Y}}f(x) {\eta^{\psi}_{\frac{1}{2}}}|_{{\bf Y}}(dx)-\int_{{\bf Y}}f(x){\eta^{\psi^n}_{\frac{1}{2}}}|_{\bf Y}(dx)\right|\nonumber\\
    &=&\left|f(1){\eta^{\psi}_{\frac{1}{2}}}|_{{\bf Y}}(\{1\})+\sum_{i=1}^\infty f(b_i){\eta^{\psi}_{\frac{1}{2}}}|_{{\bf Y}}(\{b_i\})\right.\nonumber\\
    &&\left.-\sum_{i=1}^{n}f(b_i){\eta^{\psi^n}_{\frac{1}{2}}}|_{{\bf Y}}(\{b_i\})-f(b_{n+1}){\eta^{\psi^n}_{\frac{1}{2}}}|_{{\bf Y}}(\{b_{n+1}\})\right|\nonumber\\
    &=&\left|f(1){\eta^{\psi}_{\frac{1}{2}}}|_{{\bf Y}}(\{1\})-f(b_{n+1}){\eta^{\psi^n}_{\frac{1}{2}}}|_{{\bf Y}}(\{b_{n+1}\})+\sum_{i={n+1}}^\infty f(b_i){\eta^{\psi}_{\frac{1}{2}}}|_{{\bf Y}}(\{b_i\})\right|\nonumber\\
    &\le&\left|\frac{1}{2}f(1)-\frac{1}{2}f(b_{n+1})\right|+\sum_{i={n+1}}^\infty|f(b_i)|\frac{1}{2^{i-1}}
\end{eqnarray}
where the first equality holds by (\ref{ZhangZheng2024ExampleSeparateFileEqn02}), and the second equality holds because, according to (\ref{ZhangZheng2024ExampleSeparateFileEqn02}) and (\ref{ZhangZheng2024ExampleSeparateFileEqn06}), ${\eta_{\frac{1}{2}}^{\psi^n}}|_{{\bf Y}}(\{b_i\})={\eta_{\frac{1}{2}}^\psi}|_{{\bf Y}}(\{b_i\})=\frac{1}{2^{i-1}}$ for all $i\in\{1,2,\dots,n\}$.

Since ${\bf Y}=\{b_n,~n\ge 1\}\cup\{1\}$ is endowed with the Euclidean topology, and $b_n=\sum_{i=1}^n \frac{1}{2^i}\rightarrow 1$ as $n\rightarrow \infty$, for the bounded continuous function $f$ on ${\bf Y}$, it holds that
\begin{eqnarray*}
\lim_{n\rightarrow\infty}\left|\frac{1}{2}f(1)-\frac{1}{2}f(b_{n+1})\right|=0.
\end{eqnarray*}
Moreover, since $f$ is bounded on ${\bf Y}$, we see $\lim_{n\rightarrow\infty} \sum_{i={n+1}}^\infty|f(b_i)|\frac{1}{2^{i-1}}=0$. The last two observations and (\ref{ZhangZheng2024ExampleSeparateFileEqn07}) lead to
\begin{eqnarray*}
\lim_{n\rightarrow \infty}\int_{{\bf Y}}f(x){\eta^{{\psi^n}}_{\frac{1}{2}}}|_{{\bf Y}}(dx)=\int_{{\bf Y}}f(x){\eta^{{\psi}}_{\frac{1}{2}}}|_{{\bf Y}}(dx).
\end{eqnarray*}
Since the bounded continuous function $f$ on ${\bf Y}$ was arbitrarily fixed, we see that ${\eta^{\psi^n}_{\frac{1}{2}}}|_{{\bf Y}}$ converges to ${\eta^{\psi}_{\frac{1}{2}}}|_{{\bf Y}}$ in the w-topology, as required.

Finally, we confirm that ${{\cal D}_{x_0}^\Lambda}|_{{\bf Y}} $ is not s-compact in ${\cal M}({\bf Y}).$ Suppose for contradiction that ${{\cal D}_{x_0}^\Lambda}|_{{\bf Y}} $ is   s-compact. According to \cite[Exercise 4.7.148]{Bogachev:2007}, ${{\cal D}_{x_0}^\Lambda}|_{{\bf Y}} $ is metrizable and thus sequentially s-compact. Then $\{{\eta^{\psi^n}_{\frac{1}{2}}}|_{{\bf Y}}\}_{n=3}^\infty$ has a convergent subsequence, say $\{{\eta^{\psi^{n_k}}_{\frac{1}{2}}}|_{{\bf Y}}\}_{k=1}^\infty$. Since $\{{\eta^{\psi^n}_{\frac{1}{2}}}|_{{\bf Y}}\}_{n=3}^\infty$ converges to ${\eta^{\psi}_{\frac{1}{2}}}|_{{\bf Y}}$ in the w-topology, and the s-topology is finer than the w-topology, $\{{\eta^{\psi^{n_k}}_{\frac{1}{2}}}|_{{\bf Y}}\}_{k=1}^\infty$ necessarily converges to ${\eta^{{\psi}}_{\frac{1}{2}}}|_{{\bf Y}}$ in the s-topology.
On the other hand, consider the bounded measurable function $g$ on ${\bf Y}$ defined by
\begin{eqnarray*}
    g(x)=\mathbb{I}\{x=1\}~\forall~x\in {\bf Y}.
\end{eqnarray*}
This function is not continuous on ${\bf Y}$, because ${\bf Y}=\{b_n,~n\ge 1\}\cup\{1\}$ is endowed with the Euclidean topology, and $g(b_n)\equiv 0$ does not converge to $g(1)=1$ as $n\rightarrow \infty$.  In view of (\ref{ZhangZheng2024ExampleSeparateFileEqn05}), $\int_{{\bf Y}}g(x){\eta^{\psi^{n_k}}_{\frac{1}{2}}}|_{{\bf Y}}(dx)= g(1){\eta^{\psi^{n_k}}_{\frac{1}{2}}}|_{{\bf Y}}(\{1\})= 0$ for all $k\ge 1$, and thus
\begin{eqnarray*}
   \lim_{k\rightarrow\infty}\int_{{\bf Y}}g(x){\eta^{\psi^{n_k}}_{\frac{1}{2}}}|_{{\bf Y}}(dx)=0 \neq \frac{1}{2}=g(1){\eta^{\psi}_{\frac{1}{2}}}|_{{\bf Y}}(\{1\})= \int_{{\bf Y}}g(x){\eta^{\psi}_{\frac{1}{2}}}|_{{\bf Y}}(dx).
\end{eqnarray*}
Therefore, $\{{\eta^{\psi^{n_k}}_{\frac{1}{2}}}|_{{\bf Y}}\}_{k=1}^\infty$ does not converge to ${\eta^{\psi}_{\frac{1}{2}}}|_{{\bf Y}}$ in the s-topology. This calculation incidentally proves (c), and yields a desired contradiction, justifying (d). $\hfill\Box$

\section{Conclusion}\label{ZhangZhengYear2Conclusion}
In this note, based on the recent remarkable results of Dufour and Prieto-Rumeau \cite{Dufour:2023,Dufour:2024}, we deduced that for an absorbing MDP with a given initial state, under Condition (S), the space of occupation measures has the same convergent sequences, when it is endowed with the w-topology and with the ws-topology. As a corollary, we see that under Condition (S), the w-topology and the ws-topology on the space of occupation measures coincide when the MDP is uniformly absorbing with the given initial state.  We provided two examples demonstrating that Condition (S) cannot be replaced with Condition (W), and the above assertion does not hold for the space of marginals of occupation measures on the state space. Moreover, the examples also clarify some results in the previous literature.

Finally, let us mention that there are also two natural topologies  being considered on the space $ {\cal P}_{x_0}$ of strategic measures for a general MDP model $\{{\bf X},{\bf A},p\}$, which is not necessarily absorbing: namely, the w-topology and the ws$^\infty$-topology, where the ws${}^\infty$-topology on ${\cal P}_{x_0}$ is defined as the coarsest topology, with respect to which, for each $n\ge 0$ and bounded measurable function $g$ on ${\bf H}_{n}$ that is continuous jointly in the actions $(a_1,\dots,a_{n-1})$ for each fixed $(x_0,\dots,x_n)$,   $\int_{{\bf H}_n}g(h_n)\PP|_{{\bf H}_n}(dh_n)$ is continuous in $\PP\in {\cal P}_{x_0}.$
 In \cite{Nowak:1988}, it was shown that the two topologies coincide on the space of strategic measures under Condition (S). One can check that Example \ref{ZhangZheng2024Example01} also demonstrates that this result does not hold if Condition (S) is replaced with Condition (W).

\section*{Acknowledgements} 
{We thank Alexey Piunovskiy (University of Liverpool) for useful comments, and we thank the anonymous referee for raising two interesting and relevant questions, leading to Remarks \ref{ZhangZheng2024RevisionRemark01} and \ref{ZhangZheng2024RevisionRemark02}.} 

\par\noindent\textbf{Author contributions:} Both authors contributed equally to this paper. \textbf{Yi Zhang}: Conceptualization, Investigation, Writing – original draft, Writing – review and editing. \textbf{Xinran Zheng}: Conceptualization, Investigation, Writing – original draft, Writing – review and editing.

\par\noindent\textbf{Funding sources:} This research did not receive any specific grant from funding agencies in the public, commercial, or not-for-profit sectors.
{

}


\begin{thebibliography}{00}
\bibitem{Altman:1999} Altman, E. (1999). {\em Constrained Markov Decision Processes}. Chapman and Hall/CRC, Boca Raton.
\bibitem{Balder:2001} Balder, E. (2001). On ws-convergence of product measures. {\em Math. Oper. Res.} {\bfseries 26}, 494--581.


\bibitem{Bertsekas:1978} Bertsekas, D. and Shreve, S. (1978). {\em Stochastic Optimal Control}. Academic Press, New York.

\bibitem{Bogachev:2007} Bogachev, V. (2007). {\em Measure Theory} (two volumes). Springer, Berlin.

\bibitem{Bogachev:2018} Bogachev, V. (2018). {\em Weak Convergence of Measures}. American Mathematical Society, Providence.




\bibitem{Dufour:2023} {Dufour, F. and Prieto-Rumeau, T. (2024). Absorbing Markov decision processes. {\em ESAIM: COCV}  {\bfseries 30}, article number 5.}

\bibitem{Dufour:2024} Dufour, F. and Prieto-Rumeau, T. (2024).  A note on weak compactness of occupation measures for an absorbing Markov decision process. Preprint is available at arXiv:2402.10672v1.

\bibitem{FeinbergRothblum:2012} Feinberg, E.A. and Rothblum, U. (2012). Splitting randomized stationary policies in total-reward Markov decision processes. {\em Math. Oper. Res.} \textbf{37}, 129--153.


\bibitem{FeinbergPiunovskiy:2019} Feinberg, E.A. and Piunovskiy, A. (2019). Sufficiency of deterministic policies for atomless discounted and uniformly absorbing MDPs with multiple criteria. {\em SIAM J. Control Optim.}  \textbf{57}, 163--191.

\bibitem{Ganssler:1971} G\"anssler, P. (1971). Compactness and sequential compactness in spaces of measures. {\em Z. Wahrscheinlichkeitstheorie verw Gebiete} {\bfseries 17}, 124-–146.




\bibitem{Nowak:1988} Nowak, A. (1988). On the weak topology on a space of probability measures induced by policies. {\em Bulletin Polish Acad. Sci. Math.} {\bfseries 36}, 181--186.



\bibitem{Piunovskiy:1997} Piunovskiy, A. (1997). {\em Optimal Control of Random Sequences in Problems with Constraints}, Kluwer, Dordrecht.

\bibitem{PiunovskiyZhang:2023} Piunovskiy, A. and Zhang, Y. (2024). On the continuity of  the projection mapping from strategic measures to occupation measures in absorbing Markov decision processes. {\em Appl. Math. Optim.} {\bfseries 89}, article number 58.

\bibitem{Schal:1975} Sch\"al, M. (1975). On dynamic programming: compactness of the space of policies. {\em Stoch. Proc. Appl.} {\bfseries 3}, 345--364.

\end{thebibliography}
\end{document}